\documentclass[10pt, twoside]{article}
\usepackage{amssymb,latexsym,amsmath}
\usepackage[hypertex]{hyperref}
\usepackage[english]{babel}
\date{}

\begin{document}
\makeatletter
\renewcommand{\@evenfoot}{ \thepage \hfil \footnotesize{\it ISSN 1025-6415 \ \ Reports of the National Akademy of Sciences of Ukraine, 
2002, no. 7} } 
\renewcommand{\@oddfoot}{\footnotesize{\it ISSN 1025-6415 
\ \  Dopovidi Natsionalno\" \i \ Akademi\" \i \  Nauk Ukra\"\i ni, 
2002, no. 7} \hfil \thepage } 
\par{
\leftskip=1.5cm  \rightskip=0cm  
\noindent 
UDC 512\\
\copyright \ {\bf 2002}\medskip \\
{\bf T. R. Seifullin}\medskip \\ 
{\large \bf Extension of bounded root functionals of a system\\
of polynomial equations}\bigskip \\
{\noindent \it (Presented by Corresponding Member of the NAS 
of Ukraine A. A. Letichevsky)  }  \smallskip \\
\small {\it 
The notion of a root functional of a \ system of \ polynomials or \ ideal of 
\ polynomials is a \ generalization of the notion of a root, in particular, 
for a multiple root. 
We consider bounded root functionals and 
their extension operation for a system of polynomial equation 
at which the number of equations is equal to the number of unknows.}
\par \medskip }  \noindent 
Let ${\bf R} $ be a commutative ring with unity $1$ and 
zero $0$.

 Let $x=(x_1,\ldots ,x_n)$ be variables, ${\bf R} [x]$ is the 
ring of  all  polynomials  in variables $x$ with coefficients 
in the ring ${\bf R} $.

 The degree of a monom $x^{\alpha}  = x^{ {\alpha} _1} _1\cdot 
\ldots \cdot x^{ {\alpha} _n} _n$ is called $| {\alpha} |  = 
{\alpha} _1+  \ldots   +  {\alpha} _n$, where ${\alpha} =({\alpha} 
_1,\ldots ,{\alpha} _n)$. The degree of a polynomial $F(x)$  
is  called  the  maximal degree of a monom with a nonzero coefficient, 
and such a  degree  is  denoted  by $\deg (F)$; if $F(x)=0$, then 
put $\deg (F) = -\infty $.

 {\bf Definition 1.}  Let $x=(x_1,\ldots ,x_n)$ be variables; 
we denote by ${\bf R} [x^{\leq d} ]$  the  set of all polynomials 
of degree $\leq d$. Note that ${\bf R} [x^{\leq \infty } ] = 
 {\bf R} [x]$  and  if  $d<0$,  then ${\bf R} [x^{\leq d} ] = 
\{0\}$. 

 {\bf Definition 2.}  Let $x=(x_1,\ldots ,x_n)$ be variables; 
we denote by ${\bf R} [x]_*$ the set of all maps from ${\bf R} 
[x]$ to ${\bf R} $ that are linear over ${\bf R} $,  write  such 
 maps  as  $l(x_*)$, where $x_* = (x^1_*,\ldots ,x^n_*)$, and 
 call  such  maps  linear  functionals  or  simply functionals. 
We denote the action of $l(x_*)$ on $F(x)\in {\bf R} [x]$ by 
$l(x_*).F(x).$ 

 {\bf Definition 3.}     Let    $x=(x_1,\ldots ,x_n)$     be 
    variables,     and     let $f(x)=(f_1(x),\ldots ,f_s(x))$ 
be polynomials.

 For a covector of polynomials $g(x)=(g^1(x),\ldots ,g^s(x))^{\top}$, 
we denote $f(x)g(x)= \sum\limits^{ s} _{i=1} f_i(x)g^i(x)$.

 Denote $(f(x))^{\leq d} _x = \{\sum\limits^{ s} _{i=1} f_i(x)g^i(x)| 
 \forall i=1,s: g^i(x)  \in   {\bf R} [x]$  and  
$\deg (f_i)  +\deg (g^i) \leq  d\}$.

 Denote $(f(x))_x = \{\sum\limits^{ s} _{i=1} f_i(x)g^i(x)| 
 \forall i=1,s: g^i(x) \in  {\bf R} [x]\}$.

 Note that $(f(x))^{\leq \infty } _x = (f(x))_x$, and if $d<0$, 
then $(f(x))^{\leq d} _x=\{0\}$.

 We call a functional in ${\bf R} [x]_*$ that annuls $(f(x))_x$ 
a root functional,  and a functional in ${\bf R} [x]_*$ that 
annuls $(f(x))^{\leq d} _x$ a bounded root functional. 

 {\bf Definition 4.}  Let $x=(x_1,\ldots ,x_n)$ be variables, 
and let ${\lambda} =({\lambda} _1,\ldots ,{\lambda} _n)\in {\bf 
R} ^n$;  we denote by ${\bf 1} _x({\lambda} ) = {\bf 1} _{(x_1,\ldots 
,x_n)} ({\lambda} _1,\ldots ,{\lambda} _n)$ the map such that 
${\bf 1} _x({\lambda} ).F(x) =  F({\lambda} )$ for any $F(x)\in 
{\bf R} [x]$. 

 {\bf Definition 5.}   Let  $x=(x_1,\ldots ,x_n)$,  $y=(y_1,\ldots 
,y_n)$,  and  $\hat u=(\hat u_1,\ldots ,\hat u_n)$   be variables. 
We call a difference derivative of a polynomial $F(x)\in {\bf 
R} [x]$ a covector $\hat u\nabla F(x,y) = \sum\limits^{ n} _{k=1} 
\hat u_k\nabla ^kF(x,y)$, such that
$$\hat u\nabla F(x,y) =  \sum\limits^{ n} _{k=1} 
\hat u_k\nabla ^kF(x,y) \mapsto  (x-y)\nabla F(x,y) = \sum\limits^{ 
n} _{k=1} (x_k-y_k)\nabla ^kF(x,y) = F(x) - F(y),$$
where $\forall k=1,n: \nabla ^kF(x,y) \in  {\bf R} [x,y]$. \eject

 We call a difference derivative monotonous if the degree  
of  $\nabla F(x,y)$  in $(x,y)$ is $\leq  \deg (F)-1$.

 We call a mapping that linear over ${\bf R} $ and assign, to 
a polynomial $F(x)\in {\bf R} [x]$, a covector of  a  difference 
 derivative  $\nabla F(x,y)$,  an  operator  of  difference derivative 
and denote it by $\nabla _x(x,y)$; then, we  have  $\nabla _x(x,y).F(x) 
 =  \nabla F(x,y)$. Moreover,
$$\hat u\nabla _x(x,y) = \sum\limits^{ n} _{k=1} 
\hat u_k\nabla ^k_x(x,y) \mapsto  (x-y)\nabla _x(x,y) = \sum\limits^{ 
n} _{k=1} (x_k-y_k)\nabla ^k_x(x,y) = {\bf 1} _x(x) - {\bf 1} _x(y).$$

 We call an operator of  difference  derivative  monotonous 
 if,  for  any polynomial $F(x)\in {\bf R} [x]$, the degree of 
$\nabla _x(x,y).F(x)$ is $\leq \deg (F)-1$. 

{\it  {\bf Lemma 1.}  Let $x=(x_1,\ldots ,x_n)$ and $y=(y_1,\ldots 
,y_n)$ be variables.  A  difference derivative of a polynomial 
$F(x)$ exists, for example, $\forall k=1,n:$   
$$\nabla ^kF(x,y) = \frac{  F(y_1,\ldots ,y_{k-1} ,x_k,x_{k+1} 
,\ldots ,x_n)-F(y_1,\ldots ,y_{k-1} ,y_k,x_{k+1} 
,\ldots ,x_n)} {x_k-y_k},$$
its degree is $\leq  \deg (F)-1$. A mapping that assigns, to 
 any  polynomial  $F(x)$,  a covector $\nabla F(x,y)$ is linear 
over ${\bf R} $. Thus, there exists  a  monotonous  difference 
derivative and a monotonous operator of difference derivative.

 {\bf Lemma  2.}   Let  $x=(x_1,\ldots ,x_n)$,  $y=(y_1,\ldots 
,y_n)$,   and   $\hat u=(\hat u_1,\ldots ,\hat u_n)$   be variables.

 1. For any polynomial $F(x)$, a covector $\hat u\nabla 'F(x,y) 
= \hat u\nabla F(y,x)$ is a difference derivative of the polynomial 
$F(x)$, and $\hat u\nabla '_x(x,y) = \hat u\nabla _x(y,x)$ is 
an  operator  of difference derivative.

 2. Let $V(x)=F(x)\cdot G(x)$, then $\hat u\nabla F(x,y)\cdot 
G(y)+F(x)\cdot \hat u\nabla G(x,y)$ is  a  difference derivation 
of the polynomial $V(x)$.

 3. Let $F(x) \in  {\bf R} [x^{\leq d} ]$, and let $\nabla 
'F(x,y)$ and  $\nabla ''F(x,y)$  be  two  difference derivatives 
of the polynomial $F(x)$ of degrees $\leq d-1$; then
$$\hat u\nabla 'F(x,y) = \hat u\nabla ''F(x,y) 
+ \sum\limits_{k,l} \left( (x_k-y_k)\cdot \hat u_l-(x_l-y_l)\cdot 
\hat u_k\right) \cdot T^{kl} (x,y),$$
where $k<l$ and $\deg (T^{kl} ) \leq  d-2$.} 

 {\bf Proof 1.} \begin{align*} (x-y)\cdot \nabla 'F(x,y) & = (x-y)\cdot \nabla 
F(y,x) =\\
= -(y-x)\cdot \nabla F(y,x) & = -(F(y)-F(x)) \ \ = F(x)-F(y) \end{align*}

 It follows from the first part of Statement 1 that  $\hat 
u\nabla '_x(x,y)  =  \hat u\nabla _x(y,x)$ assigns, to any polynomial 
$F(x)\in {\bf R} [x]$, its a difference derivative. The linearity 
of the map $\hat u\nabla _x(x,y)$ over ${\bf R} $ implies the linearity 
of the map $\hat u\nabla '_x(x,y) =  \hat u\nabla _x(y,x)$  over 
${\bf R} $. We finally obtain that $\hat u\nabla '_x(x,y)$ is 
an operator of difference derivative.

 {\bf Proof 2.}  \begin{align*} 
\left( (x-y)\cdot \nabla F(x,y)\right) \cdot G(y)+F(x)\cdot 
\left( (x-y)\cdot \nabla G(x,y)\right) & =\\
= \left( F(x)-F(y)\right) \cdot G(y) + F(x)\cdot 
\left( G(x)-G(y)\right) & = \\ 
=F(x)\cdot G(x)-F(y)\cdot G(y) & = V(x) - V(y).
\end{align*}

 {\bf Proof 3.}  Set $W^k(x,y) = \nabla '^kF(x,y) - \nabla 
''^kF(x,y)$, and set
$$T^{kl} (x,y) = \nabla ^k_x(x,y).W^l(x,y) = \frac{ 1} { x_k-y_k} \cdot 
(W^l(y_{<k} ,x_k,x_{>k} ,y)-W^l(y_{<k} ,y_k,x_{>k} ,y)).$$
It is directly verified that the equality in the statement 
is  true. Further, since the degrees of difference derivatives 
$\nabla 'F(x,y)$ and $\nabla ''F(x,y)$  are $\leq d-1$, then we have 
$\deg (W^l)\leq d-1$, hence, $\deg (T^{kl} ) \leq  \deg (W^l)-1 
\leq  d-2$. 

 {\bf Assumption 1.}  In the sequel, unless otherwise stated, 
we will consider only monotonous difference derivatives of polynomials 
and only monotonous operators of difference derivative.

 If $x=(x_1,\ldots ,x_n)$ are variables, then by $y\simeq x$ 
we mean $y=(y_1,\ldots ,y_n)$. \medskip

 {\it {\bf Theorem  1.}   Let  $x=(x_1,\ldots ,x_n)$  and  $y\simeq 
x$  be  variables,  let   $f(x)   =(f_1(x),\ldots ,f_n(x))$ 
 be  polynomials,  let  ${\delta} _f   =   \sum\limits^{ n} _{i=1} 
(\deg (f_i)-1)$,   and   let $F(x)\in {\bf R} [x^{\leq d} ]$; 
then we have
$$\det \left\| \begin{matrix}\nabla f(x,y) & \nabla 
F(x,y) \cr
 f(x) &  F(x)\end{matrix} \right\|  = \det \left\| \begin{matrix}\nabla 
f(x,y) & \nabla F(x,y) \cr
 f(y) &  F(y)\end{matrix} \right\| .$$
Denote this polynomial by $R(x,y)$; then we have the following:

 $R(x,y)$ have a degree $\leq {\delta} _f+d$,

 $R(x,y)$ is uniquely determined up to an addend of the form
$$\sum\limits_{i,j} \left( f_i(x)\cdot f_j(y)-f_i(y)\cdot 
f_j(x)\right) \cdot {\omega} ^{ij} (x,y)$$
indepedently of the choice of $\nabla F(x,y)$, where $i<j$ 
and  $\deg (f_i)+\deg (f_j)+\deg ({\omega} ^{ij} ) \leq  {\delta} 
_f+d$, 

$R(x,y)$ is uniquely determined up to an addend of the form
$$\sum\limits_{i,j} \left( f_i(x)\cdot f_j(y)-f_i(y)\cdot 
f_j(x)\right) \cdot {\omega} ^{ij} (x,y) + \sum\limits_i\left( 
f_i(x)\cdot F(y)-f_i(y)\cdot F(x)\right) \cdot {\Omega} ^i(x,y)$$
indepedently of the choice of $\nabla f(x,y)$, where $i<j$ 
and  $\deg (f_i)+\deg (f_j)+\deg ({\omega} ^{ij} ) 
\leq  {\delta} _f+d$ for the first summand, and $\deg (f_i) +\deg ({\Omega} ^i) 
\leq  {\delta} _f$ for the second summand. }

 {\bf Proof.}  Since $f_i(x) - \sum\limits^{ n} _{k=1} (x_k-y_k)\cdot 
\nabla ^kf_i(x,y) = f_i(x) -  (f_i(x)  -  f_i(y))  = 
f_i(y)$ and $F(x) - \sum\limits^{ n} _{k=1} (x_k-y_k)\cdot \nabla ^kF(x,y) 
= F(x) - (F(x) - F(y)) =  F(y)$, by adding, to the last 
row, the linear combination of the rest rows  of  the  first determinant 
matrix, we obtain the second determinant  matrix.  It  implies 
 the equality of determinants.

 It follows from the monotony of a  difference  derivative 
 that  the degree  of   $\nabla F(x,y)$   is   $\leq d-1$;
and the degree of $\nabla f_i(x,y)$   is   $\leq \deg (f_i)-1$ for any $i$, 
  then   the   degree   of the polynomial $R(x,y)$   is   
$\leq 
\sum\limits^{ n} _{i=1} (\deg (f_i)-1)+(\deg (F)-1) +1 \leq  {\delta} _f+d$.

 Since the degree of $\nabla F(x,y)$ is $\leq d-1$, by Statement 
3 of Lemma 2, variation of $\nabla F(x,y)$ is of the form
$$\hat u\cdot \nabla 'F(x,y) = \hat u\cdot \nabla 
F(x,y) + \sum\limits_{k,l} \left( (x_k-y_k)\cdot \hat u_l-(x_l-y_l)\cdot 
\hat u_k\right) \cdot T^{kl} (x,y),$$
where $k<l$, and $\deg (T^{kl} ) \leq  d-2$. Then $R(x,y)$ 
is uniquely determined  up  to  the addend
\begin{align*}
& \sum\limits_{k,l} \pm \det \left\| \begin{matrix}\nabla 
^{\not= k,l} f(x,y) &   0 \cr
\nabla ^kf(x,y) & -(x_l-y_l) \cr
\nabla ^lf(x,y) &  (x_k-y_k) \cr
  f(x) &   0\end{matrix} \right\| \cdot T^{kl} (x,y) = 
\hphantom{hhhhhhhhhhhhhhhhhhhhhhhhhh}
\end{align*} \begin{align*}
&\qquad =\sum\limits_{k,l} \pm \det \left\| \begin{matrix} 
          \nabla ^{\not= k,l} f(x,y) \cr
(x_k-y_k)\cdot \nabla ^kf(x,y)+(x_l-y_l)\cdot \nabla ^lf(x,y) 
\cr
f(x)\end{matrix} \right\| \cdot T^{kl} (x,y) =\\
&\qquad =\sum\limits_{k,l} \pm \det \left\| \begin{matrix}\nabla 
^{\not= k,l} f(x,y) \cr
   -f(y) \cr
    f(x)\end{matrix} \right\| \cdot T^{kl} (x,y) = \sum\limits_{i,j} 
\left( f_i(x)\cdot f_j(y)-f_i(y)\cdot f_j(x)\right) \cdot {\omega} 
^{ij} (x,y),\end{align*}
where $i<j$. The second equality is true since $\forall i=1,n:$
\begin{align*}
-f_i(y) & = (x_k-y_k)\cdot \nabla ^kf_i(x,y)+(x_l-y_l)\cdot 
\nabla ^lf_i(x,y) \ + \\
& \qquad + \sum\limits_{m\not= k,l} (x_m-y_m)\cdot 
\nabla ^mf_i(x,y) - f_i(x),\end{align*}
i. e., the last but one row of the third determinant matrix is the 
sum of the last but one row and the lineare combination of the rest row 
of the second determinant matrix. The last equality is obtained by 
decomposition of the determinant into  minors  of the two last 
rows. Moreover, we have $\deg (f_i) + \deg (f_j) + \deg ({\omega} 
^{ij} ) \leq  {\delta} _f+d$. 

 Permuting $f_t(x)$ and $F(x)$ in the statement proved above, 
 we  obtain  that $R(x,y)$ is uniquely determined up to an addend
$$\sum\limits_{i,j\not= t} \left( f_i(x)\cdot f_j(y)-f_i(y)\cdot 
f_j(x)\right) \cdot {\omega} ^{ij} (x,y) + \sum\limits_{i\not= 
t} \left( f_i(x)\cdot F(y)-f_i(y)\cdot F(x)\right) \cdot {\Omega} 
^i(x,y)$$
under lack of uniqueness of 
$\nabla f_t(x,y)$, where $i<j$ and $\deg 
(f_i) + \deg (f_j) + \deg ({\omega} ^{ij} ) \leq  {\delta} _f+d$ 
for the first summand, and $\deg (f_i) + \deg (F)  +  \deg ({\Omega} 
^i)  \leq   {\delta} _f+\deg (F)$,  hence, $\deg (f_i) + \deg 
({\Omega} ^i) \leq  {\delta} _f$, for the second  summand.  Summing 
the additional  addends appearing on changing $\nabla f_t(x$,y$)$ 
for  all  $t=1,n$,  we  obtain  that  $R(x,y)$  is uniquely 
determined up to an addend of the form
$$\sum\limits_{i,j} \left( f_i(x)\cdot f_j(y)-f_i(y)\cdot 
f_j(x)\right) \cdot {\omega} ^{ij} (x,y) + \sum\limits_i\left( 
f_i(x)\cdot F(y)-f_i(y)\cdot F(x)\right) \cdot {\Omega} ^i(x,y)$$
under lack of uniqueness of 
$\nabla f(x,y)$, where $i<j$ and $\deg (f_i) 
+ \deg (f_j) + \deg ({\omega} ^{ij} ) \leq   {\delta} _f+d$ for 
the first summand, and $\deg (f_i) +\deg ({\Omega} ^i) \leq  
{\delta} _f$ for the second summand. 
\medskip 

 {\it {\bf Theorem  2.}   Let   $x=(x_1,\ldots ,x_n)$   and   $y\simeq 
x$   be   variables,   $f(x)   =(f_1(x),\ldots ,f_n(x))$ be 
polynomials, let ${\delta} _f  =\sum\limits^{ n} _{i=1} (\deg 
(f_i)-1)$,  let  a  functional $L(x_*)$ annuls $(f(x))^{\leq 
{\delta} _f+{\delta} } _x$, where ${\delta} \geq 0$, and let 
$F(x)\in {\bf R} [x^{\leq d} ]$. We set
$$H(x) = L(y_*).\det \left\| \begin{matrix}\nabla 
f(x,y) & \nabla F(x,y) \cr
 f(x) &  F(x)\end{matrix} \right\|  = L(y_*).\det \left\| \begin{matrix}\nabla 
f(x,y) & \nabla F(x,y) \cr
 f(y) &  F(y)\end{matrix} \right\| .$$
Then we have the following:

 1. $H(x) \in  {\bf R} [x^{\leq \max ({\delta} _f,d-{\delta} 
-1)} ]$.

 2. $H(x)$ is uniquely determined up to an addend  in  $(f(x))^{\leq 
\max ({\delta} _f,d-{\delta} -1)} _x$, indepedently of the choice 
of $\nabla f(x,y)$, and uniquely determined up to an addend in 
$(f(x))^{\leq d-{\delta} -1} _x$, indepedently of the choice 
of $\nabla F(x,y)$.
\eject

 3. If $F(x)\in (f(x))^{\leq d} _x$, then $H(x)\in (f(x))^{\leq 
d-{\delta} -1} _x$.

 4.  $H(x)$  is  uniquely  determined  up  to  an  addend  
in  $(f(x))^{\leq d-{\delta} -1} _x$, indepedently of the determination 
of $L(x_*)$ outside ${\bf R} [x^{\leq {\delta} _f+{\delta} } 
]$. } 

 {\bf Proof 1.} We have
\begin{align*} 
H(x) &= L(y_*).\det \left\| \begin{matrix}\nabla 
f(x,y) & \nabla F(x,y) \cr
 f(y) &  F(y)\end{matrix} \right\|  =\\
&=L(y_*).F(y)\det \left\| \nabla f(x,y)\right\| 
 + L(y_*).\det \left\| \begin{matrix}\nabla f(x,y) & \nabla F(x,y) 
\cr
 f(y) &    0\end{matrix} \right\| .\end{align*}
The   first   summand    $\in     {\bf R} [x^{\leq {\delta} _f} ]$,    
and    the    second    summand    
$\in L(y_*).\sum\limits_{{\alpha} 
,{\beta} } (f(y))^{\leq {\alpha} } _y\cdot {\bf R} [x^{\leq {\beta} 
} ]$, where ${\alpha} +{\beta}  \leq  {\delta} _f+d$. Since  
$L(y_*)$  annuls  $(f(y))^{\leq {\delta} _f+{\delta} } _y$, without 
changing of the sum we can retain only these terms for which 
${\alpha} \geq {\delta} _f+{\delta} +1$, this means that $-{\alpha} 
\leq -({\delta} _f+{\delta} +1)$, and, hence, for the  remaining 
 terms,  we  have ${\beta} =({\alpha} +{\beta} )-{\alpha}  \leq 
 ({\delta} _f+d) - ({\delta} _f+{\delta} +1) = d-{\delta} -1$. 
Hence, the second summand $\in   \sum\limits_{\beta} {\bf R} 
[x^{\leq {\beta} } ]  \subseteq {\bf R} [x^{\leq d-{\delta} 
-1} ]$, where ${\beta} \leq d-{\delta} -1$.  Then  the  sum of 
the  both  summands  $\in   {\bf R} [x^{\leq {\delta} _f} ]  
+ {\bf R} [x^{\leq d-{\delta} -1} ] \subseteq  {\bf R} [x^{\leq 
\max ({\delta} _f,d-{\delta} -1)} ]$. Hence, we have $H(x) 
\in  {\bf R} [x^{\leq \max ({\delta} _f,d-{\delta} -1)} ]$.

 {\bf Proof 2.}  Under lack of uniqueness of $\nabla F(x,y)$, by Theorem 1,  
$H(x)$  is  uniquely determined up to an addend
$$L(y_*).\sum\limits_{i,j} f_i(x)\cdot f_j(y)\cdot 
{\omega} ^{ij} (x,y) \in  L(y_*).\sum\limits_{{\alpha} ,{\beta} 
} (f(y))^{\leq {\alpha} } _y\cdot (f(x))^{\leq {\beta} } _x,$$
where   ${\alpha} +{\beta}    \leq    {\delta} _f+d$.   The 
  last   inclusion    is    true    since    $\forall i: \deg 
(f_i)+\deg (f_j)+\deg ({\omega} ^{ij} )  \leq {\delta} _f+d$. 
 Since  $L(y_*)$  annuls  $(f(y))^{\leq {\delta} _f+{\delta} 
} _y$,   without changing the sum we can retain only these terms 
for which ${\alpha} \geq {\delta} _f+{\delta} +1$; this means 
that $-{\alpha} \leq -({\delta} _f+{\delta} +1)$, and, hence, 
for the remaining terms, we  have  ${\beta} =({\alpha} +{\beta} 
)-{\alpha}   \leq ({\delta} _f+d) - ({\delta} _f+{\delta} 
+1) = d-{\delta} -1$. Hence,  this  addend  $\in   \sum\limits_{\beta} 
(f(x))^{\leq {\beta} } _x  \subseteq (f(x))^{\leq d-{\delta} 
-1} _x$, where ${\beta} \leq d-{\delta} -1$. We finally  obtain 
 that, under lack of uniqueness of  $\nabla F(x,y)$,  $H(x)$  is uniquely 
determined up to an addend in $(f(x))^{\leq d-{\delta} -1} _x$.

 Under lack of uniqueness of   $\nabla f(x,y)$, by Theorem 1, $H(x)$ is 
uniquely determined  up to an addend of the form
$$L(y_*).\sum\limits_{i,j} f_i(x)\cdot f_j(y)\cdot 
{\omega} ^{ij} (x,y) + L(y_*).\sum\limits_i\left( f_i(x)\cdot 
F(y)-f_i(y)\cdot F(x)\right) \cdot {\Omega} ^i(x,y),$$
where $\forall i,j: \deg (f_i)+\deg (f_j)+\deg ({\omega} ^{ij} 
) \leq  {\delta} _f+d$, $\forall i: \deg (f_i) +\deg ({\Omega} 
^i) \leq  {\delta} _f$. As shown above,   the   
 first    summand    $\in     (f(x))^{\leq d-{\delta} -1} _x$. 
   Since $\sum\limits_i\left( -f_i(y)\cdot F(x)\cdot {\Omega} 
^i(x,y)\right)   \in   (f(y))^{\leq {\delta} _f} _y\cdot {\bf 
R} [x]$,  it  is  annuled  by   $L(y_*)$.   
The polynomial $\sum\limits_i\left( 
f_i(x)\cdot F(y)\cdot {\Omega} ^i(x,y)\right)  \in  {\bf R} [y]\cdot 
(f(x))^{\leq {\delta} _f} _x$. Acting  by  $L(y_*)$  on  this 
polynomial, we obtain a polynomial $\in  (f(x))^{\leq {\delta} 
_f} _x$.

 We  finally  obtain  that  this  sum \  $\in$  $(f(x))^{\leq 
d-{\delta} -1} _x   +   (f(x))^{\leq {\delta} _f} _x   \subseteq 
(f(x))^{\leq \max ({\delta} _f,d-{\delta} -1)} _x$. \ Hence, 
under lack of uniqueness of \ $\nabla f(x,y)$, \  $H(x)$  is  
uniquely determined up to an addend in 
$(f(x))^{\leq \max ({\delta} _f,d-{\delta} -1)} _x$. 
\eject

 {\bf Proof 3.}  In the proofs of Theorems 1 and 2, we  use 
 a  weaker  condition than the condition under which a difference 
derivative of the polynomial  $F(x)$ is monotonous, namely, the 
condition under which  its  degree  is  $\leq d-1$. Hence, 
these theorems are true if the last condition  is  satisfied  instead 
the first condition. 

 Let $F(x) = f(x)g(x) \in  (f(x))^{\leq d} _x$. By Statement 
2 of Lemma 2, $F(x)$ have  two difference derivatives $\nabla 
F(x,y)$ and $\nabla f(x,y)g(y)+f(x)\nabla g(x,y)$, and their 
 degrees are $\leq d-1$, although the second difference derivative 
may be not monotonous when $\deg (F)<d$ . We have \\
\begin{align*}
H(x) &= L(y_*).\det \left\| \begin{matrix}\nabla 
f(x,y) & \nabla _x(x,y).f(x)g(x) \cr
 f(y) &     f(y)g(y)\end{matrix} \right\|  \buildrel{ (f(x))^{\leq 
d-{\delta} -1} _x}\over \equiv \\
&\hbox{(by Statement 2 on the uniqueness of $H(x)$ under lack of uniqueness of
$\nabla F(x,y)$ )}\\
& \equiv L(y_*).\det \left\| \begin{matrix}\nabla f(x,y) 
& \nabla f(x,y)g(y)+f(x)\nabla g(x,y) \cr
 f(y) &    f(y)g(y)\end{matrix} \right\|  =\\
&=L(y_*).\det \left\| \begin{matrix}\nabla f(x,y) 
& f(x)\nabla g(x,y) \cr
 f(y) &     0\end{matrix} \right\|  \in  L(y_*).\sum\limits_{{\alpha} 
,{\beta} } (f(y))^{\leq {\alpha} } _y\cdot (f(x))^{\leq {\beta} 
} _x,\end{align*}
where ${\alpha} +{\beta}  \leq  {\delta} _f+d$. Since  the 
 functional  $L(y_*)$  annuls  $(f(y))^{\leq {\delta} _f+{\delta} 
} _y$,  without changing the sum, we can retain only these terms 
for which ${\alpha} \geq {\delta} _f+{\delta} +1$, this means 
that $-{\alpha} \leq -({\delta} _f+{\delta} +1)$, and, hence, 
for the remaining terms, we  have  ${\beta} =({\alpha} +{\beta} 
)-{\alpha}   \leq ({\delta} _f+d) - ({\delta} _f+{\delta} 
+1) = d-{\delta} -1$. Hence,   the  obtained  polynomial
$\in   \sum\limits_{\beta} 
(f(x))^{\leq {\beta} } _x  \subseteq (f(x))^{\leq d-{\delta} 
-1} _x$, where ${\beta} \leq d-{\delta} -1$.  Since  difference 
 of  $H(x)$  and  the  obtained  polynomial  $\in (f(x))^{\leq 
d-{\delta} -1} _x$, we have  $H(x) \in  (f(x))^{\leq d-{\delta} 
-1} _x$. 

 {\bf Proof 4.}  Let $L'(x_*) = L(x_*)$ in ${\bf R} [x^{\leq 
{\delta} _f+{\delta} } ]$, then it, as  well  as  the functional 
$L(x_*)$, annuls $(f(x))^{\leq {\delta} _f+{\delta} } _x\subseteq 
{\bf R} [x^{\leq {\delta} _f+{\delta} } ]$, and $l(x_*)  =  L'(x_*) 
 -  L(x_*)$ annuls ${\bf R} [x^{\leq {\delta} _f+{\delta} } ]$. 
We have \\
\begin{align*} &L'(y_*).\det \left\| \begin{matrix}\nabla f(x,y) 
& \nabla F(x,y) \cr
 f(x) &  F(x)\end{matrix} \right\|  - 
L(y_*).\det \left\| \begin{matrix}\nabla f(x,y) & \nabla F(x,y) \cr
 f(x) &  F(x)\end{matrix} \right\|  =\\
&\qquad =l(y_*).\det \left\| \begin{matrix}\nabla f(x,y) 
& \nabla F(x,y) \cr
 f(x) &  F(x)\end{matrix} \right\|  =\\
&\qquad =F(x)\cdot l(y_*).\det \left\| \nabla f(x,y)\right\| 
 + l(y_*).\det \left\| \begin{matrix}\nabla f(x,y) & \nabla F(x,y) 
\cr f(x) &    0\end{matrix} \right\| .\end{align*}
\\
Since $l(y_*)$ annuls ${\bf R} [y^{\leq {\delta} _f+{\delta} 
} ] \supseteq  {\bf R} [y^{\leq {\delta} _f} ]$ and $\det \left\| 
\nabla f(x,y)\right\|   \in   {\bf R} [y^{\leq {\delta} _f} ]\cdot 
{\bf R} [x]$,  the first summand is equal to $0$. The second 
 summand  $\in   l(y_*).\sum\limits_{{\alpha} ,{\beta} } {\bf 
R} [y^{\leq {\alpha} } ]\cdot (f(x))^{\leq {\beta} } _x$, where 
${\alpha} +{\beta}  \leq  {\delta} _f+d$. Since $l(y_*)$ annuls 
${\bf R} [y^{\leq {\delta} _f+{\delta} } ]$, without changing 
the sum we  can retain only these terms for which ${\alpha} \geq 
{\delta} _f+{\delta} +1$, this means that $-{\alpha} \leq -({\delta} 
_f+{\delta} +1)$, and, hence, for the remaining terms ${\beta} 
=({\alpha} +{\beta} )-{\alpha}  \leq  ({\delta} _f+d) - ({\delta} 
_f+{\delta} +1) =  d-{\delta} -1$.  Hence, the obtained polynomial 
$\in  \sum\limits_{\beta} (f(x))^{\leq {\beta} } _x \subseteq 
 (f(x))^{\leq d-{\delta} -1} _x$, where ${\beta} \leq d-{\delta} 
-1$. We  finally obtain that $H(x)$ is uniquely determined  up 
 to  an  addend  in  $(f(x))^{\leq d-{\delta} -1} _x$, indepedently 
of the determination of $L(x_*)$ outside ${\bf R} [x^{\leq {\delta} 
_f+{\delta} } ]$. 
\eject

 {\it {\bf Theorem  3.}   Let  $x=(x_1,\ldots ,x_n)$  and  $y\simeq 
x$  be  variables,  let   $f(x)   = (f_1(x),\ldots ,f_n(x))$ 
be polynomials, and  let  ${\delta} _f  =\sum\limits^{ n} _{i=1} 
(\deg (f_i)-1)$.  Let  $\forall i=1,2: L_i(x_*)$ annuls $(f(x))^{\leq 
{\delta} _f+{\delta} _i} _x$, where ${\delta} _i\geq 0$. We set
\begin{align*}
L(x_*) &= L_1(x_*).L_2(y_*).\det \left\| \begin{matrix}\nabla 
f(x,y) & \nabla _x(x,y) \cr
 f(x) &  {\bf 1} _x(x)\end{matrix} \right\|  = \\
&=L_1(x_*).L_2(y_*).\det 
\left\| \begin{matrix}\nabla f(x,y) & \nabla _x(x,y) \cr
 f(y) &  {\bf 1} _x(y)\end{matrix} \right\| .
\end{align*}
Then we have the following:

 1. $L(x_*)$ is uniquely determined in ${\bf R} [x^{\leq {\delta} 
_f+{\delta} _1+{\delta} _2+1} ]$, indepedently  of  the choice 
of $\nabla f(x,y)$ and the choice of $\nabla _x(x,y)$.

 2. $L(x_*)$ annuls $(f(x))^{\leq {\delta} _f+{\delta} _1+{\delta} 
_2+1} _x$.

 3. $L(x_*)$ is uniquely determined in ${\bf R} [x^{\leq {\delta} 
_f+{\delta} _1+{\delta} _2+1} ]$, indepedently  of  the determination 
of $L_1(x_*)$ outside ${\bf R} [x^{\leq {\delta} _f+{\delta} 
_1} ]$, and the  determination  of  $L_2(x_*)$ outside ${\bf 
R} [x^{\leq {\delta} _f+{\delta} _2} ]$. }

 {\bf Proof.}  Since $\nabla _x(x,y)$ is an operator linear 
over ${\bf R} $, $L(x_*)$ is  a  map  that linear over ${\bf 
R} $, i. e., it is a linear functional. Let a polynomial 
$F(x) \in  {\bf R} [x^{\leq {\delta} _f+{\delta} _1+{\delta} 
_2+1} ]$. Set $d={\delta} _f+{\delta} _1+{\delta} _2+1$ and ${\delta} 
 = {\delta} _2$. Then $L_2(x_*)$ annuls $(f(x))^{\leq {\delta} _f+{\delta}}_x$ and
$F(x) \in   {\bf R} [x^{\leq d} ]$. Also, we have 
$d-{\delta} -1  
=({\delta} _f+{\delta} _1+{\delta} _2+1)-{\delta} _2-1 
= {\delta} _f+{\delta} _1$,   
$\max ({\delta} _f,d-{\delta} -1)  =  
\max ({\delta} _f,{\delta} _f+{\delta} _1)  =
{\delta} _f+{\delta} _1$ since ${\delta} _1\geq 0$. Set
$$H(x) = L_2(y_*).\det \left\| \begin{matrix}\nabla 
f(x,y) & \nabla F(x,y) \cr
 f(y) &  F(y)\end{matrix} \right\| ,$$
then \\
$$L(x_*).F(x) = L_1(x_*).L_2(y_*).\det \left\| 
\begin{matrix}\nabla f(x,y) & \nabla F(x,y) \cr
 f(y) &  F(y)\end{matrix} \right\|  = L_1(x_*).H(x).$$ \\
By Statement 1 of Theorem 2, $H(x) \in  {\bf R} [x^{\leq \max 
({\delta} _f,d-{\delta} -1)} ] = {\bf R} [x^{\leq {\delta} _f+{\delta} 
_1} ]$. 

 {\bf Proof 1.}  By Statement 2 of Theorem 2, the  polynomial 
 $H(x)$  is  uniquely determined up to an addend in $(f(x))^{\leq 
d-{\delta} -1} _x = (f(x))^{\leq {\delta} _f+{\delta} _1} _x$, 
indepedently of the choice  of  $\nabla F(x,y)$,  and  is  uniquely 
 determined  up  to   an   addend   in $(f(x))^{\leq \max ({\delta} 
_f,d-{\delta} -1)} _x = (f(x))^{\leq {\delta} _f+{\delta} _1} 
_x$, indepedently of the  choice  of  $\nabla f(x,y)$. Since 
 $L_1(x_*)$  annuls $(f(x))^{\leq {\delta} _f+{\delta} _1} _x$, 
 $L(x_*).F(x)  =  L_1(x_*).H(x)$  is   uniquely determined, indepedently 
of the choice of  $\nabla _x(x,y).F(x)  =  \nabla F(x,y)$,  and 
 the choice of $\nabla f(x,y)$. From the arbitrariness of $F(x) 
\in  {\bf R} [x^{\leq {\delta} _f+{\delta} _1+{\delta} _2+1} 
]$, we obtain that $L(x_*)$ is uniquely determined  in  ${\bf 
R} [x^{\leq {\delta} _f+{\delta} _1+{\delta} _2+1} ]$,  indepedently 
 of  the choice of $\nabla f(x,y)$ and the choice of $\nabla _x(x,y)$.

 {\bf Proof 2.}  Let $F(x) \in  (f(x))^{\leq {\delta} _f+{\delta} 
_1+{\delta} _2+1} _x = (f(x))^{\leq d} _x$; then, by Statement 
3 of Theorem 2, the polynomial $H(x) \in  (f(x))^{\leq d-{\delta} 
-1} _x = (f(x))^{\leq {\delta} _f+{\delta} _1} _x$, and since 
$L_1(x_*)$ annuls $(f(x))^{\leq {\delta} _f+{\delta} _1} _x$, 
 we  have  $L(x_*).F(x)  =  L_1(x_*).H(x)  =  0$.   From   the 
arbitrariness  of  $F(x)\in (f(x))^{\leq {\delta} _f+{\delta} 
_1+{\delta} _2+1} _x$,  we  obtain   that   $L(x_*)$   annuls 
$(f(x))^{\leq {\delta} _f+{\delta} _1+{\delta} _2+1} _x$.

 {\bf Proof 3.}  By Statement 4 of Theorem 2, $H(x)$ is uniquely 
determined  up  to an addend in $(f(x))^{\leq d-{\delta} -1} 
_x = (f(x))^{\leq {\delta} _f+{\delta} _1} _x$, indepedently 
of the determination of $L_2(x_*)$ outside ${\bf R} [x^{\leq 
{\delta} _f+{\delta} } ]  =  {\bf R} [x^{\leq {\delta} _f+{\delta} 
_2} ]$.  Since  $L_1(x_*)$  annuls  $(f(x))^{\leq {\delta} _f+{\delta} 
_1} _x$, $L(x_*).F(x)  =  L_1(x_*).H(x)$  is  uniquely  determined 
  indepedently   of   the determination of $L_2(x_*)$ outside 
${\bf R} [x^{\leq {\delta} _f+{\delta} _2} ]$.

 Since the polynomial $H(x)  \in   {\bf R} [x^{\leq {\delta} 
_f+{\delta} _1} ]$,  $L(x_*).F(x)  =  L_1(x_*).H(x)$  is uniquely 
determined  indepedently  of  the  determination  of  $L_1(x_*)$ 
 outside ${\bf R} [x^{\leq {\delta} _f+{\delta} _1} ]$.

 Hence, it follows from the arbitrariness of $F(x) \in   {\bf 
R} [x^{\leq {\delta} _f+{\delta} _1+{\delta} _2+1} ]$  that the 
functional $L(x_*)$ is uniquely determined in  ${\bf R} [x^{\leq 
{\delta} _f+{\delta} _1+{\delta} _2+1} ]$,  indepedently of the 
determination of $L_1(x_*)$ outside ${\bf R} [x^{\leq {\delta} 
_f+{\delta} _1} ]$, and  the  determination  of $L_2(x_*)$ outside 
${\bf R} [x^{\leq {\delta} _f+{\delta} _2} ]$. 
\medskip

{\it {\bf Theorem  4.}   Let  $x=(x_1,\ldots ,x_n)$  and  $y\simeq 
x$  be  variables,  let   $f(x)   = (f_1(x_1,\ldots ,f_n(x))$ 
be polynomials, and  let  ${\delta} _f  =\sum\limits^{ n} _{i=1} 
(\deg (f_i)-1)$.  Let  $\forall i=1,2: L_i(x_*)$ annuls $(f(x))^{\leq 
{\delta} _f+{\delta} _i} _x$, where ${\delta} _i\geq 0$; then 
we have
$$L_1(x_*).L_2(y_*).\det \left\| \begin{matrix}\nabla 
f(x,y) & \nabla _x(x,y) \cr
 f(x) &  {\bf 1} _x(x)\end{matrix} \right\|  = L_2(x_*).L_1(y_*).\det 
\left\| \begin{matrix}\nabla f(x,y) & \nabla _x(x,y) \cr
 f(x) &  {\bf 1} _x(x)\end{matrix} \right\| $$
in ${\bf R} [x^{\leq {\delta} _f+{\delta} _1+{\delta} _2+1} 
]$. }

 {\bf Proof.}
\begin{align*}
&L_1(x_*).L_2(y_*).\det \left\| \begin{matrix}\nabla 
f(x,y) & \nabla _x(x,y) \cr
 f(x) &  {\bf 1} _x(x)\end{matrix} \right\|  = L_1(x_*).L_2(y_*).\det 
\left\| \begin{matrix}\nabla f(x,y) & \nabla _x(x,y) \cr
 f(y) &  {\bf 1} _x(y)\end{matrix} \right\|  =\\
&\hbox{(permuting $L_1(x_*)$ and $L_2(y_*)$ and substituting $x \mapsto 
 y$, $y \mapsto  x$)}\\
&\qquad =L_2(x_*).L_1(y_*).\det \left\| \begin{matrix}\nabla 
f(y,x) & \nabla _x(y,x) \cr f(x) &  {\bf 1} _x(x)\end{matrix} \right\|= \\
&\qquad =L_2(x_*).L_1(y_*).\det 
\left\| \begin{matrix}\nabla 'f(x,y) & \nabla '_x(x,y) \cr
  f(x) &  {\bf 1} _x(x)\end{matrix} \right\|  =\\
&\qquad =L_2(x_*).L_1(y_*).\det \left\| \begin{matrix}\nabla 
f(x,y) & \nabla _x(x,y) \cr f(x) &  {\bf 1} _x(x)\end{matrix} \right\|.
\end{align*}
By Statement 1 of Lemma 2, $\nabla '_x(x,y) =  \nabla _x(y,x)$ 
 is  operator  of  a  difference derivative, $\forall i=1,n: 
\nabla 'f_i(x,y) = \nabla f_i(y,x)$ is a  difference  derivative 
 of  the polynomial $f_i(x)$. The last   equality  is  
true  in  ${\bf R} [x^{\leq {\delta} _f+{\delta} _2+{\delta} 
_1+1} ]$  by Statement 1 of Theorem 3, since this  functional 
 is  uniquely  determined  in ${\bf R} [x^{\leq {\delta} _f+{\delta} 
_1+{\delta} _2+1} ]$, indepedently of the choice  of  $\nabla 
_x(x,y)$  and  the  choice  of $\nabla f(x,y)$. Hence, the both 
functionals is coincide in ${\bf R} [x^{\leq {\delta} _f+{\delta} 
_1+{\delta} _2+1} ]$. \bigskip 

{\footnotesize

\begin{enumerate}

\item {\it Seifullin, T. R.} 
Root functionals and root polynomials 
of a system of polynomials. (Russian)
Dopov. Nats. Akad. Nauk Ukra\"\i ni  -- 1995, -- no. 5, 5--8.

\item {\it Seifullin, T. R.} Root functionals and root relations 
of a system of polynomials. (Russian) 
Dopov. Nats. Akad. Nauk Ukra\"\i ni  -- 1995, -- no. 6, 7--10.

\item  {\it Seifullin, T. R.}  Homology of the Koszul complex of a 
system of polynomial equations. (Russian)
Dopov. Nats. Akad. Nauk Ukr. Mat. Prirodozn. Tekh. Nauki 1997, no. 9, 43--49. 

\item {\it Seifullin, T. R.} Koszul complexes of embedded systems of 
polynomials and duality. (Russian) 
Dopov. Nats. Akad. Nauk Ukr. Mat. Prirodozn. 
Tekh. Nauki 2000, no. 6, 26--34.  

\item  {\it Seifullin, T. R.}  Koszul complexes of systems of 
polynomials connected by linear dependence. (Russian) 
Some problems in 
contemporary mathematics (Russian), 326--349, Pr. Inst. Mat. Nats. Akad. Nauk 
Ukr. Mat. Zastos., 25, Natsional. Akad. Nauk Ukra\"\i ni, Inst. Mat., Kiev, 
1998.
\\
\end{enumerate}

\small{\noindent
{\it V. M. Glushkov Institute of Cybernetics of the NAS of Ukraine, Kiev
\hfill Received 06.07.2001
} 
\\ \\
}

\end{document}